\documentclass[runningheads,a4paper]{llncs}

\usepackage{amssymb}
\usepackage{graphicx}

\newcommand{\eref}[1]{(\ref{#1})}

\newcommand{\by}{\mathbf{y}}

\newcommand{\q}{\mathbf{q}}

\newcommand{\tKi}{\tilde K^{(i)}}
\newcommand{\Ki}{ K^{(i)}}
\newcommand{\tKj}{\tilde K^{(j)}}
\newcommand{\Kj}{ K^{(j)}}
\newcommand{\tP}{\tilde P}
\newcommand{\tPi}{\tilde P^{(i)}}
\newcommand{\Pii}{ P^{(i)}}
\newcommand{\ta}{\tilde a}
\newcommand{\txi}{\tilde \xi}
\newcommand{\txii}{\tilde \xi^{(i)}}
\newcommand{\xii}{ \xi^{(i)}}
\newcommand{\tomega}{\tilde \omega}

\begin{document}

\mainmatter  

\title{Numerical Methods for the Optimal Control of Scalar Conservation Laws}

\titlerunning{Numerical Methods for the Optimal Control of Scalar Conservation Laws}

%
%
\author{Sonja Steffensen \and Michael Herty%
\thanks{RWTH Aachen University, Templergraben 55, D-52065 Aachen, GERMANY.
{\tt \{herty,steffensen\}@mathc.rwth-aachen.de}}%
\and  Lorenzo Pareschi
\thanks{University of Ferrara, Department of
Mathematics, Via Machiavelli 35, I-44121 Ferrara, ITALY. {\tt
lorenzo.pareschi@unife.it}}}
\authorrunning{S. Steffensen. M. Herty and L. Pareschi}

\institute{RWTH Aachen University, Templergraben 55, \\D-52065 Aachen, GERMANY\\
{\tt \{herty,steffensen\}@mathc.rwth-aachen.de}\\
University of Ferrara, Department of
Mathematics,\\ Via Machiavelli 35, I-44121 Ferrara, ITALY\\
{\tt lorenzo.pareschi@unife.it}
}

%
%

\toctitle{Lecture Notes in Computer Science}
\tocauthor{Authors' Instructions}
\maketitle

\begin{abstract}
We are interested in a class of numerical schemes for the optimization 
of nonlinear hyperbolic partial differential equations. We present  continuous
and discretized relaxation schemes for scalar, one-- conservation
laws. We present  numerical 
results on  tracking type problems with nonsmooth desired states
and convergence results for higher--order spatial and temporal
discretization schemes. 
\end{abstract}

\section{Introduction}
We consider an optimal control problem for scalar conservation laws of the 
type 
\begin{equation}\label{ocp-orig}
 \begin{array}{rl}
{\rm minimize}_{ u_0} & J(u(T),u_0) \\[0.5em]
{{\rm subject ~to}}&
u_t +f(u)_x = 0, \qquad u(0,x) = u_0(x),
\end{array}
\end{equation}
Here, $J$ and $f$ are assumed to be smooth and possibly
nonlinear functions. The initial value $u_0$ acts as control
to the problem.  It can be observed that the wave 
interactions that occur in the solution $u$ in the case of a nonlinear 
flux function $f$ pose  the serious analytical challenges. Recently,
the differentiability of $J$ with respect to $u_0$ could be
proven in the sense of shift--differentiability. We refer to 
\cite{Bianchini2000,BressanGuerra1997,BressanLewicka1999,BressanMarson1995,Ulbrich2001,Ulbrich2002d,Ulbrich2003,BandaHerty2010,Zuazua2008} for more details.
\par
Here, a class of numerical methods applied to the optimal control 
problem (\ref{ocp-orig}) is studied. We only consider the  case of smooth 
initial data and smooth solutions $u$ and refer to 
\cite{BandaHerty2010} for more details. For a numerical analysis including
shock  waves and in the case of the Lax--Friedrichs scheme we
refer to \cite{Ulbrich2009,Zuazua2008} and the references therein. 

\subsection{Relaxation Method}

As motiviation for a numerical scheme we follow the ideas of Jin and Xin \cite{JinXin1995}. 
Therein, a  linear approximation  (\ref{imex}) 
of the nonlinear hyperbolic equation 
\[\partial_{t}u+ \partial_{x} f(u) =0\]
has been discussed. 
For   initial conditions $u(x,0)=u_{0}$ the approximation is 
\begin{equation}\label{imex}
\begin{array}{rcll}
 \partial_{t} u + \partial_{x}v&=&0\,,&u(x,0)=u_{0},\\[0.5em]
\partial_{t} v + a^2 \partial_{x}u &=& \frac{1}{\epsilon} \left( f(u)-v \right),\quad &v(x,0) = f(u_{0})\\[0.5em]
\end{array}
\end{equation}
where $\varepsilon>0$ is the relaxation rate and 
$a$ is a given constant satisfying the subcharacteristic condition
$ \max_{u} |f'(u)|\leq a.$
 For $\varepsilon$ being small, the solution $u$ of (\ref{imex}) 
satisfies 
$\partial_{t} u + \partial_{x} f(u)=
\varepsilon\partial_x((a^2-f(u)^2)\partial_{x} u)$ (cf. \cite{JinXin1995}).
Applying the relaxation to the optimal control problem \eref{ocp-orig}, we obtain 
 \begin{equation}\label{ocp-imex}
    \min_{u_{0}} J( u(\cdot,T),u_0 ) \quad \mbox{ subject to }
  \left\{\begin{array}{l}
         u_t + v_x =0,\\[0.3em]
         v_t + a^{2} u_x = \frac{1}{\epsilon} \left(
        f(u) -v \right),\\[0.3em]
        u(0,x)=u_{0},\,\,v(0,x) =
        f(u_{0})
         \end{array}\right.
 \end{equation}
The corresponding adjoint equations for \eref{ocp-imex} are given by (cf. \cite{BandaHerty10})
\begin{eqnarray*}
 - p_t - a^2 q_x = \frac{q}\epsilon  f'(u), && \; p(T,x)
 =p_{T}(x), \\
 -q_t -  p_x  = -\frac{q}\epsilon , && \; q(T,x) = q_{T}(x).
 \end{eqnarray*}
For more information on the relaxation system, its limiting scheme for $\epsilon=0$, 
further numerical analysis and extensions we refer to \cite{MR1766856,NatAreg96,BandaSeaid2001,GottliebShuTadmor2001,ARS97,BPR11,JinXin1995,MR1855278,PR} 
and the references therein.  Also, the computations are valid provided that
all appearing functions are at least once differentiable. This is in general
not the case for conservation laws.

\section{IMEX-Runge-Kutta Discretization}

Numerical discretization of the relaxation system using higher order temporal discretizations
combined with  higher order spatial discretization has been investigated in several recent
publications as for example \cite{JinXin1995,PR}. We apply so called  implicit--explicit Runge-Kutta methods \cite{PR03,PR,ARS97}
as temporal discretization (IMEX RK). Here, the expliciti integration is used for the linear hyperbolic transport part and an implicit method is applied to the  
the stiff source term. Implicit-explicit Runge-Kutta method have been studied in the context
of control problems for example in \cite{BandaHerty2010,HertySchleper11}.  Define
\[\by=(u,v)^T, \qquad g(\by)=(v,a^2
u)^T \qquad \mbox{and}\qquad  r(\by):=(0,-(v-f(u)))^T
\]
then \eref{imex} becomes
\[
\by_t + g(\by)_x = \frac{1}{\varepsilon} r(\by),  \qquad \mbox{   and  
} \qquad \by(0,x) = (u^0,f(u^0))^T(x)
\]
Applying a suitable discretization $D_x$  of the spatial derivative
yields the semi-discrete state equations
\begin{equation}\label{y-ode}
\by' = - D_x\,g(\by) + \frac{1}{\epsilon}r(\by) ,
 \qquad \by(0) =\by^0.
\end{equation}

\begin{remark} 
Spatial discretizations for the linear transport part are well--known. The simplest possible is a first--order
Upwind method:
\[
\frac{\partial}{\partial t} \by_j =-  \frac{1}{\Delta x} 
\left(\begin{array}{cc}
      0 &1\\a^2 &0
\end{array}\right) (\by_{j+1/2}-\by_{j-1/2}) + \frac{1}{\epsilon}r(\by_j)\,,
\]
where $\by_{j+1/2}$ is obtained by applying the first-order upwind method to characteristic
 variables $v \pm a u$.  Higher order MUSCL schemes,  WENO schemes or central schemes
 have also been studied in this context.
\end{remark}
The resulting semi--discrete optimal control problem is then given by:
\begin{equation}\label{ocp-1}
\begin{array}{rl}
{\rm minimize }\quad & j(\by(T),\by^0) \\[1.5em]
{{\rm subject ~to } }\quad &
\by' =  - D_x\,g(\by) + \frac{1}{\epsilon}r(\by) ,
 \qquad \by(0) =\by^0.\quad t\in[0,T]
\end{array}
\end{equation}
In the context of relaxation schemes the  semi--discrete problem is seen
as a time--integration problem with stiff source which is discretized by an IMEX RK methods. For the numerical discretization
we therefore consider the previous problem as an optimal control problem 
 involving ordinary differential equations. Literature concerning the numerical analysis of 
  Runge-Kutta methods  for the optimality system of (\ref{ocp-1}) have been studied in  \cite{Hager00,BonnansLaurent-Varin06,LangVerwer2011}. In \cite{BonnansLaurent-Varin06,Hager00}
  partitioned 
Runge-Kutta methods for the  optimality system are obtained using the {\em discretize--then--optimize} 
approach. The derived partitioned Runge--Kutta methods have been analysed with regard
to symplecticity and order of convergence. 
In \cite{HertySchleper11}, Herty and Schleper, moreover, analysed the associated adjoint
 imex Runge-Kutta method that one obtains if an explicit method is applied to  $D_x g(y)$ 
and  a (diagonally) implicit method to  $\frac{1}{\epsilon}r(y)$. In the following, we will 
analyse general partitioned Runge-Kutta methods using IMEX RK methds. More details
can be found in \cite{HertyPareschiSteffensen}.  Therein, 
the following IMEX Runge-Kutta discretization of \eref{y-ode} is studied.
\begin{equation}\label{SRK1}
\begin{array}{rcl}
Y^{(i)}_n&  =&y_n + h\sum_{j=1}^{i-1} \ta_{ij}D_x g(Y_n^{(j)}) +
h\sum_{j=1}^i a_{ij}\frac{1}{\epsilon}r(Y_n^{(j)})\qquad
i=1,..,s\\[0.6em] 
y_{n+1}&=& y_n +h\sum_{i=1}^s \tomega_iD_x
g(Y_n^{(i)}) + h\sum_{i=1}^s
\omega_i\frac{1}{\epsilon}r(Y_n^{(i)}),\qquad n=0,1,2,.
 \end{array} 
\end{equation}
A nonlinear variable transformation
and two intermediate states $\tKi_n$ and $\Ki_n$
give the equivalent system
\begin{equation}\label{SRK2}
 \begin{array}{rcl}
\tKi_n&  =& D_x g\left(y_n+h\sum_{j=1}^{s} \ta_{ij}\tKj_n +
h\sum_{j=1}^s a_{ij}\Kj_n\right)\qquad i=1,..,s\\[0.3em]
 \Ki_n&  =&
\frac{1}{\epsilon}r\left(y_n+h\sum_{j=1}^{s}  \ta_{ij}\tKj_n +
h\sum_{j=1}^s a_{ij}\Kj_n\right)\qquad i=1,..,s\\[0.3em] 
y_{n+1}&=& y_n
+h\sum_{i=1}^s \tomega_i\tKi_n + h\sum_{i=1}^s
\omega_i\Ki_n,\qquad n=0,1,2,.
 \end{array} 
\end{equation}
The associated optimality systems for the two previous optimization problems then
coincide and we refer to \cite{HertyPareschiSteffensen} for mor details. It is proven that the 
 adjoint schemes are equivalent to  
 \begin{eqnarray}
\tPi&  =&p_n - h\sum_{j=1}^s \tilde \alpha_{ij}\,g'(Y^{(j)}_n)^T\bar
D_x \tP^{(j)} - h\sum_{j=1}^s
\alpha_{ij}\,\frac{1}{\epsilon}r'(Y^{(j)}_n)^T P^{(j)}\qquad i=1,..,s \nonumber \\[0.5em] 
\Pii&  =&p_n-h\sum_{j=1}^s \tilde
\beta_{ij}\,g'(Y^{(j)}_n)^T\bar D_x\tP^{(j)} - h\sum_{j=1}^s
\beta_{ij}\, \frac{1}{\epsilon}r'(Y^{(j)}_n)^T P^{(j)}\qquad
i=1,..,s \label{ARK} \\[0.5em]
 p_{n+1}&=& p_{n} - h\sum_{i=1}^s \tomega_i g'(Y^{(i)}_n)^T\bar D_x
\tPi - h\sum_{i=1}^s \omega_i \frac{1}{\epsilon}r'(Y^{(i)}_n)^T
\Pii\qquad n=0,1,..,N-1 \nonumber
 \end{eqnarray}
Here, the coefficients of the  Runge-Kutta method
 $\tilde \alpha_{ij},\alpha_{ij}, \tilde \beta_{ij}$ and $\beta_{ij}$
are given by
\[
 \tilde \alpha_{ij}:= \tomega_j-\frac{\tomega_j}{\tomega_i} \ta_{ji}, \qquad
\alpha_{ij}:= \omega_j-\frac{\omega_j}{\tomega_i} \ta_{ji},\qquad 
\tilde \beta_{ij}:= \tomega_j-\frac{\tomega_j}{\omega_i} a_{ji}, \qquad 
\beta_{ij}:= \omega_j-\frac{\omega_j}{\omega_i} a_{ji} .\qquad \qquad
\]

\subsection{Properties of Discrete IMEX-RK Optimality System}
For the resulting scheme (\ref{SRK1}),(\ref{ARK}) order conditions can 
be stated \cite{HertyPareschiSteffensen}. To this end we add  a suitable equation
for $\tilde p$ to the previous system.
\begin{equation}\label{tilde p}
 \tilde p_{n+1}= \tilde p_{n} - h\sum_{i=1}^s \tomega_i \,f_y(Y^{(i)}_n)^T \tPi
-h \sum_{i=1}^s \omega_i \,g_y(Y^{(i)}_n)^T \Pii. 
\end{equation}
The full method therefore is a standard additive Runge-Kutta scheme for 
 \begin{eqnarray*}
 \by' &=& - D_x\,g(\by) + \frac{1}{\epsilon}r(\by)\\
\tilde{\bf p}'&=&g'(\by)^TD_x\tilde{\bf p} + \frac{1}{\epsilon}r'(\by)^T\bf p\\
\bf p'&=&g'(\by)^TD_x\tilde{\bf p} + \frac{1}{\epsilon}r'(\by)^T \bf p
\end{eqnarray*}
 If we  define 
\begin{eqnarray*}
 c_i&:=\sum_{j=1}^s a_{ij}, \qquad \qquad \mbox{and} \qquad \qquad
\tilde c_i&:=\sum_{j=1}^s \ta_{ij}, \\
\gamma_i& := \sum_{j=1}^s \alpha_{ij}, \qquad \qquad \mbox{and} \qquad \qquad
 \tilde \gamma_i &:= \sum_{j=1}^s \tilde\alpha_{ij},\\
\delta_i& := \sum_{j=1}^s \beta_{ij}, \qquad  \qquad \mbox{and} \qquad \qquad
 \tilde \delta_i &:= \sum_{j=1}^s \tilde\beta_{ij}
\end{eqnarray*}
then (\ref{thm1}) holds true. 
\begin{theorem} \label{thm1}
Consider the Runge-Kutta  scheme \eref{SRK1},\eref{ARK},\eref{tilde p}. This scheme is of 
\begin{itemize}
\item \textbf{First-Order} : if (SRK1) is of first order\\[1em]
\item   \textbf{Second-Order} : if (SRK1) is of second order\\[1em]
 \item  \textbf{Third-Order} : if (SRK1) is of third order and either
\[
\sum_{i=1}^s \omega_i \,\gamma_i ^2=\frac{1}{3},\qquad 
\sum_{i=1}^s \omega_i \,\tilde \gamma_i ^2=\frac{1}{3},\qquad
\sum_{i=1}^s \omega_i \,\gamma_i \tilde \gamma_i=\frac{1}{3},
\]
are satisfied or if 
\[
\sum_{i=1}^s \omega_i \,a_{ij}\,\gamma_i =\frac{1}{6},\qquad 
\sum_{i=1}^s \omega_i \, \ta_{ij}\, \tilde \gamma_i=\frac{1}{6}
\]
and if
\[
\sum_{i=1}^s \omega_i  \,a_{ij}\,\tilde \gamma_i =\frac{1}{6}\qquad \mbox{or} \qquad 
\sum_{i=1}^s \omega_i  \,\ta_{ij}\,\gamma_i =\frac{1}{6}
\]
\end{itemize}
are satisfied.
\end{theorem}
Note that  the system \eref{SRK1} and \eref{ARK}
is not completely coupled, since the forward scheme \eref{SRK1} is solved independently of the
adjoint scheme \eref{ARK}. General order conditions can be found e.g. in \cite{KC03}. 
The proof of Theorem \ref{thm1} and together with more details
are discussed in \cite{HertyPareschiSteffensen}.

\section{Numerical Results}

\subsection{Scalar Example}
As a simple example, we use a tracking type functional  $J(u)$ together with Burgers' equation
\[
u_t+\left(\frac{u^2}{2}\right)_x=0,
\] 
and the desired state $u_d$ at final time $T=2.0$, that belongs to the initial 
condition $u_d(0,x)=\frac{1}{2}+\sin(x)$ and we start the optimization with the initial 
data $u^s(0,x)\equiv 0.5$. Moreover, the spatial interval is given
by $x\in [0,2\pi],$
As discretization of the objective functional, we use $$J(u(\cdot,T),u_0, u_{d}) = \frac{\Delta x}{2}\sum\limits_{i=1}^K
\|u_i - u_{d\,i}\|^2.$$ Moreover, the discrete gradient of the reduced cost functional is given by 
$$\nabla_{u_{0,i}} \tilde J = p_{0,i} + (Df(u_0)^{T} \q_{0})_i.$$ In order to solve the optimal control problem, we
apply a steepest descent method (with respect to the reduced cost functional) with fixed stepsize $0<\alpha<1$,
i.e. we set $u_0^{k+1}=u_0^k+ \alpha \nabla_{u_{0,i}} \tilde J$. As stopping criterion for the optimization process
we test $|\tilde J(u_0,u_d) |<tol$ where $tol=1E-2$ denotes a predefined stopping tolerance. We observe
grid independence in the case where $u$ and $u_0$ are differentiable in space and time. \\
As first-order scheme, we test the Implicit-Explicit Euler scheme
 \[\begin{array}{rcl}
  u_i^*&=& u_i^n \\[0.2em]
  v_i^*&=&v_i^n-\frac{\Delta t}{\epsilon}(v_i^*-f(u^*_i))\\[0.2em]  
  u_i^{n+1}&=&u_i^*-\Delta t\, D_x v_i^* \\[0.2em]  
  v_i^{n+1}&=&v_i^*-\Delta t \,a^2 D_x v_i^* 
   \end{array}
\]
for the forward, as well as for the backward
\[\begin{array}{rcl}
q_i^*&=&q_i^{n+1} -\Delta t D_x^* p_i^{n+1}\\[0.2em]
p_i^*&=&p_i^{n+1}-\Delta t \,a^2 D_x^* q_i^{n+1}\\[0.2em]
q_i^{n}&=&q_i^*-\frac{\Delta t}{\epsilon} q_i^n\\[0.2em]
p_i^{n}&=&p_i^* + \frac{\Delta t}{\epsilon}q_i^n f'(u^n_i) 
\end{array}
\]
The  spatial gridsize is chosen to be $N_x=300,$ whereas the 
 time discretization is done according to the CFL condition with 
constant $c_{CFL}=0.5$.

\vspace{1em}    
\begin{center}
\includegraphics[height=0.35\textwidth]{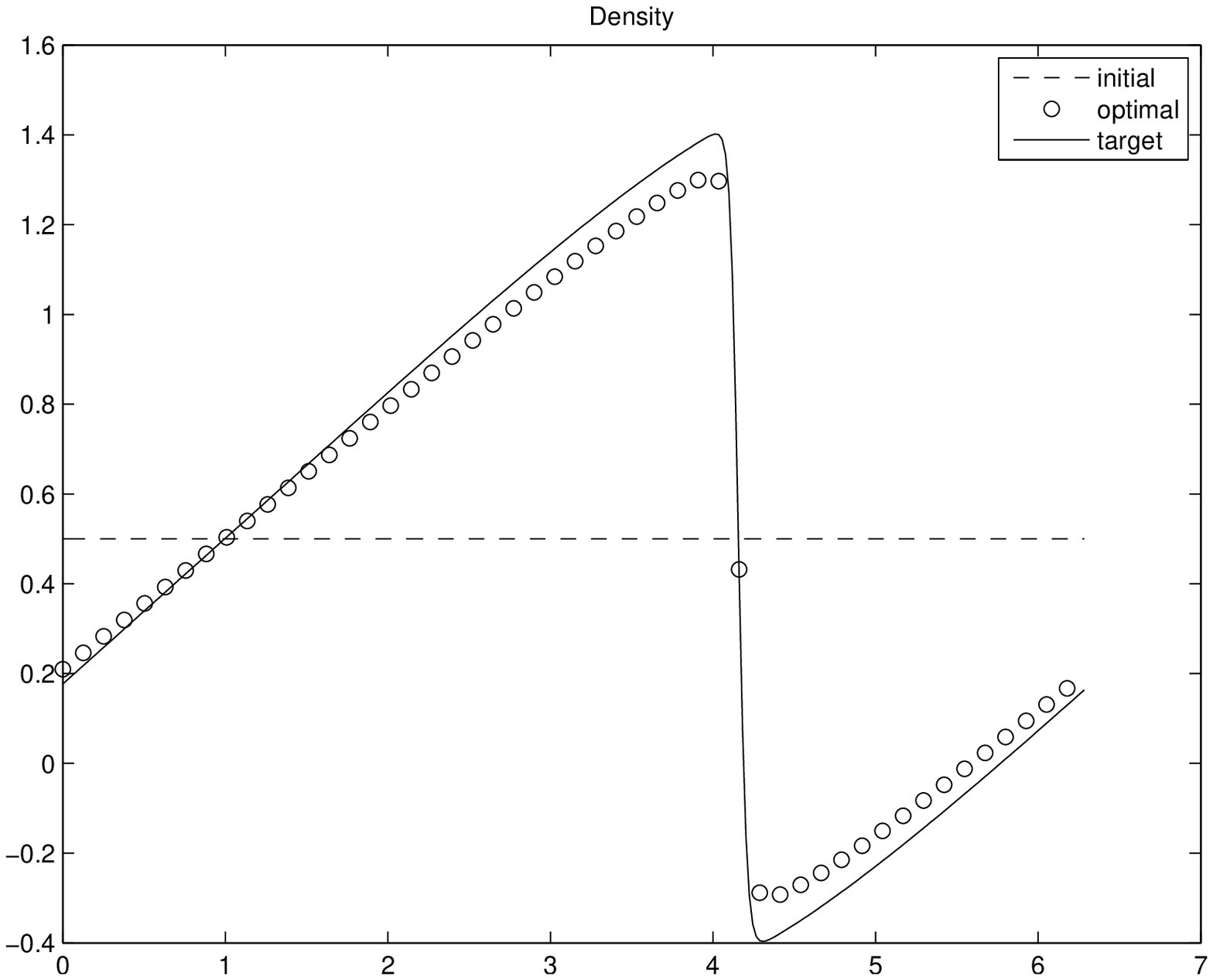} 
\end{center}  

\begin{center}
\begin{tabular}{|c|c|c|}
\hline
$\quad$ N$\quad$ & Nr. of It. & $\quad$ CPU time (in sec.) $\quad$\\
\hline\hline
100  & 44  &  1.795572e+01\\ 
150  & 43  &  3.768984e+01\\
200  &  42  &  6.380441e+01 \\ 
300  &  41  &  1.491838e+02 \\
  \hline  
\end{tabular}
\end{center}

\section{Summary}
We briefly discussed a class of numerical methods applied to an optimal control
problem for scalar, hyperbolic partial differential equations. Order conditions
for the temporal numerical discretization in the case of differentiable functions
have been stated. Future work includes the analysis of additional properties
of the derived numerical discretizations as for example strong stability 
and asymptotic preservation properties. 

\subsection*{Acknowledgments.} 
This work has been supported by DFG
HE5386/7-1,  HE5386/8-1 and by DAAD  50727872, 50756459 and
54365630. We also acknowledge the support of Ateneo Italo-Tedesco
(AIT) under the Vigoni project 2010-2012 'Adjoint implicit-
-explicit methods for the numerical solution to optimization
problems'.

\section*{Appendix.}
The discrete adjoint equations that correspond to the discrete optimization problem
associated with \eref{SRK2} are
\begin{eqnarray*}
\txii_n&  =&h\,\tomega_i \,p_{n+1}+h\sum_{j=1}^s
\ta_{ji}g'(Y^{(j)}_n)^T\bar D_x \txi^{(j)}_n + h\sum_{j=1}^s
\ta_{ji}\frac{1}{\epsilon}r'(Y^{(j)}_n)^T   
\xi^{(j)}_n\\[0.5em] 
\xii_n&  =&h\,\omega_i \,p_{n+1}+h\sum_{j=1}^s
a_{ji}g'(Y^{(j)}_n)^T\bar D_x\txi^{(j)}_n + h\sum_{j=1}^s
a_{ji}\frac{1}{\epsilon}r'(Y^{(j)}_n)^T \xi^{(j)}_n\\
[0.5em] p_{n}&=&
p_{n+1} +\sum_{i=1}^s  g'(Y^{(i)}_n)^T\bar D_x \txii _n +\sum_{i=1}^s 
\frac{1}{\epsilon}r'(Y^{(i)}_n)^T \cdot\xii_n\\ p_N&=&j'(y_N,y^0)\,.
 \end{eqnarray*} 
Moreover, the variable transformation that is needed to obtain \eref{ARK} is given by 
\[
 \tPi_n:= \frac{\txii_n}{h \,\tomega_i} \qquad \mbox{and} \qquad  \Pii_n:= \frac{\xii_n}{h \,\omega_i} 
\qquad (i=1,..,s; \quad n=0,..,N-1)\,.
\]
On the other hand, using \eref{SRK1} the associated discrete adjoint equations are 
 \begin{eqnarray*}
\zeta^{(i)}_n&  =&h\left(\tomega_i f_y(Y_n^{(i)})+ \omega_i g_y(Y_n^{(i)})\right)^Tp_{n+1}
+ \sum_{j=1}^{s} \ta_{ji}\,f_y(Y_n^{(i)}) ^T\zeta^{(j)}\\[0.3em]
&&\qquad + h\sum_{j=1}^s
a_{ji}\,g_y(Y_n^{(i)})^T \zeta^{(j)} \qquad \qquad i=1,..,s\nonumber \\[0.5em]
p_{n}&=& p_{n+1} + \sum_{i=1}^s\zeta^{(i)}_n, \qquad  \qquad \; i=1,..,N-1, \qquad   p_N=j'(y_N)
 \end{eqnarray*}
which can be transformed into the scheme \eref{ARK} using the variable transformation
\[
 \tPi_n:= p_{n+1}+\sum_{j=1}^s \frac{\ta_{ji}}{\tomega_i}\zeta^{(j)}_n\qquad \mbox{and} \qquad
\Pii_n:=p_{n+1}+\sum_{j=1}^s \frac{a_{ji}}{\omega_i}\zeta^{(j)}_n\,.
\]

.

\end{document}